\documentclass{article}
\usepackage[named]{algo}
\usepackage{spconf}
\usepackage{graphicx}
\usepackage{amsmath}
\usepackage{amsfonts}
\usepackage{amssymb}
\usepackage{epsfig}
\usepackage{url}

\newcommand{\bo}[1]{{\bf #1}}

\newcommand{\bPsi}{\mbox{\boldmath $\Psi$}}

\newcommand{\bm}[1]{\mbox{\boldmath $#1$}}

\newcommand{\boldeta}{\mbox{\boldmath  $\eta$}}

\title{Deconvolution of Poissonian Images Using \\Variable
Splitting and Augmented Lagrangian Optimization}

\name{ M\'{a}rio A. T. Figueiredo \hspace{1.5cm} Jos\'{e} M. Bioucas-Dias }
\address{{\it Instituto de Telecomunica\c{c}\~{o}es}, \\
{\it Instituto Superior T\'{e}cnico}, Technical University of Lisbon, {\bf Portugal}\\
Email: $\{$mario.figueiredo,\  jose.bioucas$\}$@lx.it.pt}

\begin{document}

\maketitle

\begin{abstract}
 Although much research
has been devoted to the problem of restoring Poissonian images,
namely in the fields of medical and astronomical imaging, applying
the state of the art regularizers (such as those based on wavelets
or total variation) to this class of images is still an
open research front. This paper proposes a new image deconvolution approach for
images with Poisson statistical models, with the following building blocks:
(a) a standard regularization/MAP criterion, combining the Poisson
log-likelihood with a regularizer (log-prior) is adopted; (b) the
resulting optimization problem (which is difficult, since it involves
a non-quadratic and non-separable term plus a non-smooth term) is
transformed into an equivalent constrained problem, via a variable
splitting procedure; (c) this constrained problem is addressed using
an augmented Lagrangian framework. The effectiveness of the resulting
algorithm is illustrated in comparison with current state-of-the-art
methods.
\end{abstract}

\section{Introduction}
\subsection{Poissonian Images}
Image restoration is one of the earliest and most classical
inverse problems in imaging, dating back to the 1960's. Much of the work
in this field has been devoted to developing regularizers (priors
or image models, in a Bayesian perspective) to deal with the ill-conditioning
or ill-posedness of the observation operator, and to devising
efficient algorithms to solve the resulting optimization problems.

A large fraction of the work on image restoration assumes that
the observation operator is linear  (often the convolution
with some blur point spread function) and the presence of additive
Gaussian noise. For this scenario, recent work has lead to a set
of state-of-the-art restoration methods, which involve non-smooth
convex regularizers ({\it e.g.}, total-variation, $\ell_1$ norm of frame coefficients)
and efficient special-purpose algorithms (see \cite{TwIST}, \cite{CombettesSIAM},
\cite{FigueiredoNowak2003}, \cite{Starck2003b}, and references therein).

The algorithms developed for the linear/Gaussian observation model
cannot be directly applied to other statistical ({\it e.g.}, Poisson
or Gamma) observation models. The Poisson case is
well studied and highly relevant in fields such as astronomical \cite{StarckBook},
biomedical \cite{Dey}, \cite{FesslerHero},  and photographic
imaging \cite{Foi_ICIP2005}. A very recent overview of deconvolution
methods for Poissonian images can be found in \cite{DupeFadiliStarck}, where
a state-of-the-art algorithm is also introduced.

Although our approach can be applied to other regularizers, we focus
here on total-variation (TV), well-known for its discontinuity
preserving ability \cite{Chambolle04}, \cite{art:Rdin:O:F:Physica:92}.
The combination of TV regularization with the log-likelihood resulting
from the Poissonian observations of a convolved image, leads to
an objective function with a non-quadratic non-separable term
(the log-likelihood) plus a non-smooth term (TV).
This objective function poses the following
difficulties to the current state-of-the-art algorithms: (a) the Poisson log-likelihood term
doesn't have a Lipschitz-continuous gradient, which
is a necessary condition for the applicability of algorithms of the
forward-backward splitting (FBS) class \cite{CombettesSIAM}, \cite{DupeFadiliStarck};
(b) the presence of a  convolution in the observation model precludes the
direct application of the Douglas-Rachford splitting methods
described in \cite{CombettesPesquet}. Moreover, if an FBS algorithm is
applied (ignoring that the convergence conditions are not met), it
is known to be slow, specially when the observation operator is
ill-conditioned, a fact which has stimulated recent
research aimed at obtaining faster methods \cite{FISTA}, \cite{TwIST},
\cite{SpaRSA_SP}.

In this paper, we propose a new approach to tackle the optimization
problem referred to in the previous paragraph. Firstly, the original
optimization problem  is transformed into an equivalent constrained one,
via a variable splitting procedure. Secondly, this constrained problem
is addressed using an algorithm developed within the augmented Lagrangian
framework, for which convergence is guaranteed. The effectiveness of the
resulting algorithm is illustrated in comparison with current
state-of-the-art alternatives \cite{DupeFadiliStarck}, \cite{Foi_ICIP2005},
\cite{Dey}.

\section{Augmented Lagrangian}
In this section,  we briefly review the augmented Lagrangian framework,
a key building block of our approach.
Consider a convex optimization problem with linear equality constraints
\begin{equation}\begin{array}{cl}
 {\displaystyle \min_{{\bf v}\in \mathbb{R}^d}} & E({\bf v})\\
 \mbox{s.t.} & {\bf A v} = {\bf  b}, \end{array}\label{constrained_linear}
\end{equation}
where ${\bf b} \in \mathbb{R}^p$ and ${\bf A}\in \mathbb{R}^{p\times
d}$. The so-called augmented Lagrangian function for this problem is defined
as
\begin{equation}
{\cal L}_A ({\bf v},\boldeta,\mu) = E({\bf v}) + \boldeta^T ({\bf
Av-b}) + \frac{\mu}{2}\,  \|{\bf Av-b}\|_2^2,\label{augmented_L}
\end{equation}
where $\boldeta \in \mathbb{R}^p$ is a vector of Lagrange
multipliers and $\mu \geq 0$ is called the AL penalty
parameter \cite{Nocedal}. The AL algorithm iterates
between minimizing ${\cal L}_A ({\bf v},\boldeta,\mu)$ with respect
to ${\bf v}$, while keeping $\boldeta$ fixed, and updating $\boldeta$.

\vspace{0.1cm}
\begin{algorithm}{AL}{}
Set $k=0$, choose $\mu > 0$, ${\bf v}_0$, and  $\boldeta_0$.\\
\qrepeat\\
     ${\bf v}_{k+1} \in \arg\min_{{\bf v}} {\cal L}_A ({\bf v},\boldeta_k,\mu)$\\
     $\boldeta_{k+1} \leftarrow \boldeta_{k} + \mu ({\bf Av}_{k+1} - {\bf b})$\\
     $k \leftarrow k + 1$
\quntil stopping criterion is satisfied.
\end{algorithm}
\vspace{0.1cm}

It is possible (in some cases recommended) to update the
value of $\mu$ at each iteration \cite{Nocedal}.
Notice, however, that it is not necessary to take $\mu$ to
infinity to guarantee convergence to the solution of the
constrained problem (\ref{constrained_linear}).
In this paper, we will consider only the case of fixed $\mu$.

After a straightforward manipulation, the terms
added to $E({\bf v})$ in  ${\cal L}_A ({\bf v},\boldeta_k,\mu)$
(see (\ref{augmented_L})) can be written as a single quadratic term,
leading to the following alternative form for the AL algorithm:
\begin{algorithm}{AL (version 2)}{}
Set $k=0$, choose $\mu > 0$, ${\bf v}_0$, and  ${\bf d}_0$.\\
\qrepeat\\
     ${\bf v}_{k+1} \in \arg\min_{{\bf v}} E({\bf v}) + \frac{\mu}{2}\|{\bf Av-d}_k\|_2^2$\\
     ${\bf d}_{k+1} \leftarrow {\bf d}_{k} - ({\bf Av}_{k+1} - {\bf b})$\\
     $k \leftarrow k+1$
\quntil stopping criterion is satisfied.
\end{algorithm}
\vspace{0.1cm}

This form of the AL algorithm makes clear its equivalence with the
recently introduced Bregman iterative method \cite{YinOsherGoldfarbDarbon}.

\section{Problem Formulation}
Let $\bo{y} = (y_1,...,y_n) \in \mathbb{N}_{0}^{n}$ denote an $n$-elements observed image or signal of
counts, assumed to be a sample of a random image $\bo{ Y} = (Y_1,...,Y_n) \in \mathbb{N}_0^n$
composed of $n$ independent Poisson variables
\begin{equation}
P[\bo{Y} = \bo{y} | \bm{\lambda} ] = \prod_{i=1}^n \frac{\lambda_i^{y_i} \, e^{-\lambda_i}}{y_i !},
\label{Poisson}
\end{equation}
where $\bm{\lambda} = (\lambda_1,...,\lambda_n) \in \mathbb{R}_+^n$ is the underlying mean signal,
assumed to be a blurred version of an unknown ${\bf x}$, {\it i.e.},
\begin{equation}
\bm{\lambda} = {\bf K\, x},\label{blur}
\end{equation}
where ${\bf K}$ is the matrix representation of the blur operator,
which is herein assumed to be a convolution.
When dealing with images, we adopt the usual vector notation obtained
by stacking the pixels into an $n$-vector using,
e.g., lexicographic order. Combining (\ref{Poisson}) and (\ref{blur}), we can write
\[
\log P[\bo{Y} = \bo{y} | {\bf x} ] = \sum_{i=1}^n y_i\,
\log\left( ({\bf K\, x})_i\right)-({\bf K\, x})_i - \log(y_i!)
\]
where $({\bf K\, x})_i$ denotes the $i$-th component of ${\bf K\, x}$ \cite{Dey}, \cite{StarckBook}.

Under the regularization or the Bayesian maximum a posteriori (MAP) criterion,
the original image ${\bf x}$ is inferred by
solving a minimization problem with the form
\begin{eqnarray}
       \min_{\bo{x}} & &  L(\bo{x})  \label{eq:L_uncons}\\
       \mbox{s.t.} & & \bo{x} \geq {\bf 0}.\label{eq:L_uncons_2}
\end{eqnarray}
The objective  $L(\bo{x})$ is the penalized negative log-likelihood,
\begin{eqnarray}
L(\bo{x}) & = &  -\log P[\bo{Y} = \bo{y} | {\bf x} ] + \tau\, \phi(\bo{x}) ,\label{criterion1}\\
& = & \sum_{i=1}^n ({\bf K\, x})_i - y_i\, \log\left( ({\bf K\, x})_i\right)  + \tau\, \phi(\bo{x}),
\label{criterion2}
\end{eqnarray}
where  $\phi:\mathbb{R}^n\rightarrow \mathbb{R}$ is the penalty/regularizer (negative of the log-prior, from  the Bayesian perspective),
and $\tau \in \mathbb{R}_+$ is the regularization parameter. Notice that the non-negativity constraint on ${\bf x}$
guarantees that $\bm{\lambda} = {\bf K\, x}$  is also non-negative, if all the entries in
${\bf K}$ are non-negative (as is the case in most convolution kernels modeling a variety of
blur mechanisms).

In this  work, we adopt the TV regularizer \cite{Chambolle04}, \cite{art:Rdin:O:F:Physica:92}, {\it i.e.},
\begin{equation}
 \phi(\bo{x})  =  \mbox{TV}(\bo{x}) = \sum_{s=1}^n \sqrt{(\Delta^h_s\bo{x})^2+(\Delta^v_s\bo{x})^2},\label{eq:theTV}
\end{equation}
where  $(\Delta^h_s\bo{x}$ and $\Delta^v_s\bo{x})$ denote the horizontal and vertical first
order differences at pixel $s\in\{1,\dots,n\}$, respectively.

Each term $({\bf K\, x})_i - y_i\, \log\left( ({\bf K\, x})_i\right)$
of (\ref{criterion2}),  corresponding to the negative log-likelihood,
is  convex, thus so is their sum. If the space of constant images $\{ {\bf x} = \alpha(1,1,...,1),\;
 \alpha\in\mathbb{R}\}$,
for which TV is zero, does not belong to the null space of ${\bf K}$, and
the counts  $(y_1,...,y_n)$ are all non-zero, then the objective function $L$ is
coercive and strictly convex thus possessing a unique minimizer
\cite{CombettesSIAM}.

\section{Proposed Approach}
\subsection{Variable Splitting}
The core of our approach consists in rewriting the optimization problem
defined by (\ref{eq:L_uncons})--(\ref{eq:theTV})
as the following equivalent constrained problem:
\begin{eqnarray}
\min_{{\bf x,z,u}} & & \sum_{i=1}^n (z_i - y_i\, \log z_i) + \tau \phi({\bf u})\\
\mbox{s. t.} & & {\bf K\, x} = {\bf z}\\
& & {\bf x} = {\bf u}.
\end{eqnarray}
Notice that we have dropped the non-negativity constraint (\ref{eq:L_uncons_2});
this constraint could be applied to either ${\bf x}$, ${\bf z}$, or ${\bf u}$ (as long
as all elements of ${\bf K}$ are non-negative). However, as shown below, if
applied to ${\bf z}$, this constraint will be automatically satisfied during the execution
of the algorithm, thus can be dropped. Notice that this problem can be written compactly
in the form (\ref{constrained_linear}),
using the translation table
\begin{equation}
{\bf v} = \left[\begin{array}{c} {\bf x}\\ {\bf z}\\ {\bf u} \end{array}\right], \hspace{0.5cm}
{\bf b} = {\bf 0}, \hspace{0.5cm}  {\bf A} = \left[\begin{array}{ccc} {\bf K} & -{\bf I} & {\bf 0}\\
{\bf I} & {\bf 0} & -{\bf I}
 \end{array}\right],
\end{equation}
and with
\begin{equation}
E({\bf v}) = E({\bf x,z,u}) = \sum_{i=1}^n (z_i - y_i\, \log z_i) + \tau \phi({\bf u}).
\end{equation}

\subsection{Applying the AL Algorithm}
The application of Step 3 of the AL (version 2) algorithm to the problem
just described requires the solution of a joint minimization with respect to ${\bf x}$,
${\bf z}$, and ${\bf u}$, which is still a non-trivial problem.
Observing that each partial minimization (e.g., with respect to ${\bf x}$,
while keeping ${\bf z}$ and ${\bf u}$ fixed) is computationally treatable suggests
that this joint minimization can be addressed using a non-linear block
Gauss-Seidel (NLBGS) iterative scheme. Of course, this raises the
question of wether such a scheme converges, and of how much computational
effort ({\it i.e.}, iterations) should be spent in solving
this minimization in each step of the AL algorithm. Experimental evidence
({\it e.g.} \cite{GoldsteinOsher}) suggests that good results are
obtained by running just one NLBGS step in each step
of the AL algorithm. In fact, it has been shown that the AL algorithm
with a single NLBGS step per iteration does converge
\cite{EcksteinBertsekas}, \cite{Setzer}.
Remarkably, the only condition required is that the objective function
be proper and convex.

Finally, applying AL  (version 2), with a single NLBGS step per iteration,
to the constrained problem presented in the previous subsection leads to our
proposed algorithm, termed PIDAL ({\it Poisson image deconvolution by AL}).
The algorithm is presented in Fig.~\ref{fig:PIDAL}.

\begin{figure}
\begin{center}
\fbox{\parbox{0.95\columnwidth}{\small
\begin{algorithm}{Poisson Image Deconvolution by AL (PIDAL)}{}
Choose $\bo{x}_0$, $\bo{z}_0$, $\bo{u}_0$, $\bo{d}^{(1)}_{0}$, $\bo{d}^{(2)}_{0}$,
$\mu$, and $\tau$. Set $k := 0$.\\
\qrepeat \\
  ${\bf x}' = {\bf z}_k+\bo{d}^{(1)}_{k}$\\
  ${\bf x}'' = {\bf u}_k+\bo{d}^{(2)}_{k}$\\
  ${\displaystyle \bo{x}_{k+1} := \arg\min_{\bo{x}} \| {\bf K\, x} -{\bf x}'\|_2^2} +
                   \| {\bf x} -{\bf x}''\|_2^2$ \\
        ${\bf z}' = {\bf Kx}_{k+1}-\bo{d}^{(1)}_{k}$\\
${\displaystyle \bo{z}_{k+1} := \arg\min_{\bo{z}} \sum_{i=1}^n z_i - y_i\, \log z_i
   + \frac{\mu}{2}\| {\bf z} -{\bf z}'\|_2^2 }$\\
   ${\bf u}' = {\bf x}_{k+1}-\bo{d}^{(2)}_{k}$\\
${\displaystyle\bo{u}_{k+1} := \arg\min_{\bo{x}} \frac{1}{2}\|
\bo{u}-\bo{u}'\|^2+ (\tau/\mu)\,\phi(\bo{u})}$.\\
$\bo{d}^{(1)}_{k+1}  :=  \bo{d}^{(1)}_{k} - (\bo{K\, x}_{k+1} -\bo{z}_{k+1})$\\
$\bo{d}^{(2)}_{k+1}   :=  \bo{d}^{(2)}_{k} - (\bo{x}_{k+1} -\bo{u}_{k+1})$\\
$k := k + 1$
\quntil some stopping criterion is satisfied.
 \end{algorithm} }}
 \caption{The PIDAL algorithm.}\label{fig:PIDAL}
\end{center}
\end{figure}

The minimization with respect to $\bo{z}$ (line 5) is given by
\begin{equation}\label{eq:wiener_fillter}
    {\bf x}_{k+1} = \left({\bf K}^T{\bf K}+{\bf I}\right)^{-1}({\bf K}^T{\bf x}'+{\bf x}'').
\end{equation}
We are assuming that $\bf K$ models a convolution, thus it is a block Toeplitz or
block circulant matrix and (\ref{eq:wiener_fillter}) can be implemented in $O(n\log n)$ operations,
using the FFT algorithm.

Step 7 is separable and has closed form: for each  $z_i$, it amounts to computing
the non-negative root of the second order polynomial $\mu z_i^2+(1-\mu\, z_i')z_i - y_i$, given by
\begin{equation}\label{eq:second_order_roots}
      z_{i,k+1} = \left(\mu\, z_i'-1+\left((\mu\, z_i'-1)^2+4\,\mu\, y_i\right)^{1/2}\right)/(2\mu).
\end{equation}
Notice that this is always a non-negative value, thus justifying the statement
made above that the constraint ${\bf z}\geq 0$ is automatically satisfied by
the algorithm.

The minimization with respect to $\bo{u}$ (line 9) is, by definition,
the Moreau proximity mapping $\bPsi_{\tau\phi}:\mathbb{R}^n\rightarrow \mathbb{R}^n $
of the regularizer $\tau \phi$ \cite{CombettesSIAM}.
In this paper, the adopted regularizer is the TV norm (\ref{eq:theTV}), thus
${\bf u}_{k+1}$ is obtained by applying TV-based denoising to ${\bf u}'$.
To implement this denoising operation, we use Chambolle's well-known
algorithm \cite{Chambolle04}, although other fast methods
are also available \cite{WangYangYinZhang}.

Notice how the variable splitting, followed by the augmented Lagragian
approach, converted a difficult problem (\ref{eq:L_uncons})--(\ref{eq:theTV}),
involving a non-quadratic and non-separable term plus a (non-smooth) TV regularizer,
into a sequence of three simpler problems: (a) quadratic problem with a linear solution (line 5);
(b) a separable problem with closed-form solution (line 7); (c) a
TV-based denoising problem (line 9), for which efficient algorithms exist.

\section{Experiments}
We now report experiments where  PIDAL is compared with two
state-of-the-art methods \cite{DupeFadiliStarck}, \cite{Foi_ICIP2005}.
All the experiments use synthetic data produced according to
(\ref{Poisson})--(\ref{blur}), where ${\bf x}$ is the {\em Cameraman}
image and ${\bf K}$ represents a uniform blur.
In Experiment 1 (following \cite{Foi_ICIP2005}), the blur is $9\times 9$,
and the original image is scaled to a maximum value of 17600; this is a
high SNR situation. In Experiment 2 (following \cite{DupeFadiliStarck})
the blur is $7\times 7$, and the maximum value of ${\bf x}$ belongs to
$\{5,30, 100, 255\}$; this represents low SNR situations.

Parameter $\mu$ of the PIDAL algorithm affects its convergence speed,
but its adaptive choice is a topic beyond the scope of this paper. In all
the experiments, we use $\mu = \tau/50$, found to be a good rule of thumb.
PIDAL is initialized
with ${\bf x}_0 = {\bf y}$, ${\bf z}_0 = {\bf K\, x}_0$, ${\bf u}_0 = {\bf x}_0$,
${\bf d}^{(1)}_0 = {\bf 0}$, and ${\bf d}^{(2)}_0 = {\bf 0}$.

In Experiment 1, the regularization parameter $\tau$ was set to $6\times 10^{-4}$;
since our goal is to propose a new algorithm, not a new deconvolution
criterion, we didn't spend time fine tuning $\tau$ or using methods to
adaptively estimate it from the data. Since the method in \cite{Foi_ICIP2005}
includes a set of adjustable parameters which need to be hand tuned, the
comparison remains fair. The improvement in SNR (ISNR) obtained by PIDAL was 6.96dB (average
over 10 runs), better than the 6.61dB reported in \cite{Foi_ICIP2005}. This
result is more remarkable if we notice that the TV regularizer is considerably
simpler than the locally adaptive  approximation
techniques used in \cite{Foi_ICIP2005}.

\begin{table}
\centering
\caption{Mean absolute errors obtained by PIDAL and the algorithm from \cite{DupeFadiliStarck} (average
over 10 runs).}\label{tab:results}
\vspace{0.3cm}
\begin{tabular}{|l  || c | c | c | c |}
\hline
max intensity & 5 & 30 & 100 & 255 \\ \hline \hline
PIDAL & 0.37 & 1.34 & 3.99 & 8.65 \\ \hline
Algorithm from \cite{DupeFadiliStarck} & 0.44 & 1.44 & 4.69 & 10.40 \\
\hline
\end{tabular}
\end{table}

For Experiment 2, we downloaded the code available at {\tt\small www.greyc.ensicaen.fr/$\sim$fdupe/.}
Although the regularizer is not the same, we used the same values of $\tau$ found
in that code; if anything, this constitutes a disadvantage
for PIDAL. Following \cite{DupeFadiliStarck}, the accuracy of
an image estimate $\widehat{\bf x}$ is assessed by the mean absolute error MAE $= \|\widehat{\bf x}-{\bf x}\|_1/n$.
Table~\ref{tab:results} shows the MAE values achieved by PIDAL and the algorithm
of \cite{DupeFadiliStarck}, for the several values of the maximum original image
intensity, showing that PIDAL always yields lower MAE. In our experiments, each
run of the algorithm from \cite{DupeFadiliStarck} takes roughly 10 times longer than PIDAL.

\section{Concluding Remarks}
We have proposed an approach to TV deconvolution
of Poissonian images, by exploiting a variable splitting
procedure and augmented Lagrangian optimization.
In the experiments reported in the paper, the proposed algorithm
exhibited state-of-the-art performance. We are currently working
on extending our methods to other regularizers, such as
those based on frame-based  sparse representations.

\end{document}